\documentclass[11pt]{article}

\usepackage{latexsym}

\begin{document}

%%%%%%%%%%%%%%%%%%%%%%%%%%%%%%%%%%%%%%%%%%%%%%%%%%%%%%%%%%%% begin lama.sty
%  LAMA.STY  (auszug)                                           6. M\"arz 92
%

% Spezielle mathematische Symbole ----------------------------------------
\newcommand{\EE}{\mathop{\rm I\! E}\nolimits}
\newcommand{\E}{\mathop{\rm E}\nolimits}
\newcommand{\I}{\mathop{\rm Im  }\nolimits}
\newcommand{\R}{\mathop{\rm Re  }\nolimits}
\newcommand{\CC}{\mathop{\rm C\!\!\! I}\nolimits}
\newcommand{\FF}{\mathop{\rm I\! F}\nolimits}
\newcommand{\KK}{\mathop{\rm I\! K}\nolimits}
\newcommand{\LL}{\mathop{\rm I\! L}\nolimits}
\newcommand{\MM}{\mathop{\rm I\! M}\nolimits}
\newcommand{\NN}{\mathop{\rm I\! N}\nolimits}
\newcommand{\PP}{\mathop{\rm I\! P}\nolimits}
\newcommand{\QQ}{\mathop{\rm I\! Q}\nolimits}
\newcommand{\RR}{\mathop{\rm I\! R}\nolimits}
\newcommand{\ZZ}{\mathop{\rm Z\!\!Z}\nolimits}
% ------------------------------------------------------------------------
\newcommand{\integer}{\mathop{\rm int}\nolimits}
\newcommand{\erf}{\mathop{\rm erf}\nolimits}
\newcommand{\diag}{\mathop{\rm diag}\nolimits}
\newcommand{\fl}{\mathop{\rm fl}\nolimits}
\newcommand{\eps}{\mathop{\rm eps}\nolimits}
\newcommand{\var}{\mathop{\rm var}\nolimits}

%%%%%%%%%%%%%%%%%%%%%%%%%%%%%%%%%%%%%%%%%%%%%%%%%%%%%%%%%%%% end lama.sty

\newcommand{\bild}[3]{{        % \bild{H"ohe}{Unterschrift}{ref-label}
  \unitlength1mm
  \begin{figure}[ht]
  \begin{picture}(120,#1)\end{picture}
  \caption{\label{#3}#2}
  \end{figure}
}}

\newcommand{\ring}{{\cal K}}

\newcommand{\pfeil}{\rightarrow}

\newcommand{\WD}{{\rm WD}}
\newcommand{\WB}{{\rm WB}}
\newcommand{\ZD}{{\rm ZD}}
\newcommand{\HD}{{\rm HD}}
\newcommand{\BD}{{\rm BD}}
\newcommand{\HA}{{\rm HA}}
\newcommand{\HB}{{\rm HB}}
\newcommand{\BA}{{\rm BA}}
\newcommand{\ZA}{{\rm ZA}}
\newcommand{\ZB}{{\rm ZB}}
\newcommand{\HBB}{{\rm\cal HB}}
\newcommand{\BAB}{{\rm\cal BA}}
\newcommand{\BBfull}{{\rm BB}}
\newcommand{\BB}{{\rm BB}}
\newcommand{\BBB}{{\rm\cal BB}}
\newcommand{\BBs}{{\rm {\widetilde{BB}}}}
\newcommand{\AK}{{\rm AK}}
\newcommand{\tr}{{\rm tr}}
\newcommand{\Tr}{{\rm Tr}}
\newcommand{\iso}{\stackrel{\sim}{=}}

\newcommand{\kat}{{\cal C}}
\newcommand{\rmat}{{\cal R}}
\newcommand{\oalg}{{\cal A}}
\newcommand{\falg}{{\cal F}}
\newcommand{\eich}{{\cal G}}
\newcommand{\hilb}{{\cal H}}
\newcommand{\calm}{{\cal M}}
\newcommand{\mod}{{\cal M}}

\newcommand{\bigrho}{\rho_\oplus}
\newcommand{\bigphi}{\phi_\oplus}

\newtheorem{hypo}{Hypothese}
\newenvironment{bew}{Proof:}{\hfill$\Box$}
\newtheorem{bem}{Remark}
\newtheorem{bsp}{Example}
\newtheorem{axiom}{Axiom}
\newtheorem{de}{Definition}
\newtheorem{satz}{Proposition}
\newtheorem{lemma}[satz]{Lemma}
\newtheorem{kor}[satz]{Corollary}
\newtheorem{theo}[satz]{Theorem}
\newtheorem{satzdef}[satz]{Definition and Proposition}
 
\newcommand{\sbegin}[1]{\small\begin{#1}}
\newcommand{\send}[1]{\end{#1}\normalsize}
\newcommand{\myownlab}[1]{\label{#1}}

\sloppy
%%%%%%%%%%%%%%%%%%%%%%%%%%%%%%%%%%%%%%%%%%%%%%%%%%%%%%%%%%%%%%%%%%%%%%%%%

\newcommand{\lang}{1}

\title{The Birman-Murakami-Wenzl Algebra of Coxeter Type D}
\author{Reinhard Oldenburg\\
email: reinhard.oldenburg@math.uni-augsburg.de}
\date{October 10, 2016}
\maketitle

\abstract{The paper defines a generic Birman-Wenzl algebra of Coxeter Type D and investigates its structure as a semi-simple algebra.}

\begin{de}
The braid group $\ZD_n$ of Coxeter D-type has generators $X_0,\ldots,X_{n-1}$
with relations
\begin{eqnarray}
X_iX_jX_i&=&X_jX_iX_j\qquad i,j\geq1,|i-j|=1\label{hd1}\\
X_0X_2X_0&=&X_2X_0X_2\label{hd2}\\
X_0X_j&=&X_jX_0\qquad j\neq2\label{hd40} \\
X_iX_j&=&X_jX_i\qquad |i-j|>1,i,j\geq1\label{hd4} 
\end{eqnarray}
In addition, the Coxeter Group $\WD_n$ generators satisfy the quadratic relation  for all $i\geq0$
\begin{eqnarray}
X_i^2=1
\end{eqnarray}
\end{de}

\begin{de}
The braid group $\ZB_n$ of Coxeter B-type has 
generators $Y,X_1,\ldots,X_{n-1}$ and relations
\begin{eqnarray}
YX_1YX_1&=&X_1YX_1Y\label{hb1}\\
X_iX_jX_i&=&X_jX_iX_j\qquad |i-j|=1\label{hb2} \\
YX_i&=&X_iY\qquad i\geq2\label{hb3} \\
X_iX_j&=&X_jX_i\qquad |i-j|>1  \label{hb4} 
\end{eqnarray}
In addition, the Coxeter Group $\WB_n$ includes quadratic relations:
\begin{eqnarray}
X_i^2&=&1								  \\ Y^2&=&1
\end{eqnarray}
\end{de}

\begin{lemma}\label{morph}
There is an injective morphism of Coxeter  groups $i:\WD_n\rightarrow\WB_n$ which
maps $X_0\mapsto YX_1Y,X_i\mapsto X_i$.
\end{lemma}
\begin{bew}
The proof that $i$ is a morphism is simple.
We only check (\ref{hd3}):
\[i(X_0X_1)=YX_1YX_1=X_1YX_1Y=i(X_1)i(X_0)=i(X_1X_0)\]
and (\ref{hd2})
\begin{eqnarray*} i(X_0X_2X_0)&=&YX_1YX_2YX_1Y=YX_1X_2Y^2X_1Y=
YX_1X_2X_1Y\\	
&=&YX_2X_1X_2Y=	X_2YX_1YX_2=i(X_2X_0X_2)
\end{eqnarray*}
As all relations preserve the parity of the number of $Y$ in the words, it is clear that none of the words with an odd number of $Y$ are in the image. On the other hand, words in $\WB_n$ with an even number of $Y$ lie in the image: Chose  a word in $\WB_n$ that has an even number of $Y$. Then we can choose the first and second occurrence of $Y$ and using (\ref{hb3}) we can bring them together until only a power of $X_1$ is in between so that we have the sub-word $w=YX_1^kY$. If $k$ is even, then $w=1$, else $w=YX_1Y=i(X_0)$. Hence, the subgroup of $\WB_n$ consisting of words with an even number of $Y$ is isomorphic to $\WD_n$.
\end{bew}

\begin{de}		
Let $q$ denote an invertible element in an 
integral domain $R$. The D-type Hecke algebra $\HD_n$
over this ring is the quotient of the group algebra 
of $\ZD_n$ by the relation
\begin{eqnarray}
X_i^2&=& (q-1)X_i+q, i\geq 0\label{hd} 
\end{eqnarray}
\end{de}

\begin{de}			
Let $q,p_0,p_1$ denote  invertible elements in an 
integral domain $R$. The B-type Hecke algebra $\HB_n$
over this ring is the quotient of the group algebra of $\ZB_n$ by the relations
\begin{eqnarray}
X_i^2&=& (q-1)X_i+q\label{hb5} \\
(Y-p_0)(Y-p_1)&=& 0 \label{hb6}
\end{eqnarray}
\end{de}

Note that (\ref{hb6}) is equivalent to $Y^2=(Q-1)Y+Q$ with  $p_0=Q,p_1=-1$ and this is a more common way of giving the quadratic relation.    
				   
We would like to extend the morphism from lemma \ref{morph} to the Hecke algebras.

\begin{lemma}
The map $X_0\mapsto kYX_1Y,X_i\mapsto X_i$ for some $k\in R$
defines an injective morphism $i:\HD_n\rightarrow\HB_n$ if $p_0+p_1=0$ and $k={{-1} \over {p_0p_1}}$.
\end{lemma}
\begin{bew}
Under the assumption on the parameters stated in the lemma the quadratic relation for $Y$ simplifies to $Y^2=-p_0p_1=1/k$.
We have  to check that $i(X_0)=kYX_1Y$ satisfies the correct quadratic relation:
$(kYX_1Y)^2=k^2YX_1Y^2X_1Y=kYX_1^2Y=(q-1)kYX_1Y+qkY^2=(q-1)kYX_1Y+q$.
Bases of Hecke algebras may be labeled by words in the Coxeter
group. Now, $i$ is an injection in the group case and, due
to the relation $Y^2=1/k$, the image of  $i$ closes as a
subalgebra of $\HD_n$. It consists of those words that have
an even number of $Y$. 
\end{bew}

By the results of Ariki and Koike \cite{ariki},\cite{ariki1}
the algebra $\HB_n$ is semi-simple in the generic case. 
The simple $\HB_n$ Modules are parametrized by the set of
pairs $T=(T_0,T_1)$ of Young diagrams of total size $n$.
For such a pair $T$ the simple module $V_T$ has a basis of
pairs of Young tableaux $t=(t_0,t_1)$ of shape $T$.
The generator $Y$ acts as multiplication with its eigenvalues (i.e. if $(Y-p_0)(Y-p_1)=0$, then as multiplication with $p_0$ resp. $p_1$) depending on the tableaux containing 1, i.e. if $1$ is contained
in $t_i$ it acts as multiplication with  $p_i$.
$X_i$ acts as multiplication with $q$ if $i$ and $i+1$ 
are contained in the same row of the same tableaux, and as $-1$
if they are contained in the same column of the same tableaux.
Otherwise the result of $X_i$ acting on   $t$
is a linear combination of $t$ and $t'$, where in $t'$ 
the numbers $i,i+1$ are interchanged.

Using the morphism $i$ we can view $V_T$  as a module of $\HD_n$ given a restriction on the eigenvalues of $Y$. 
\begin{lemma}
Suppose that the eigenvalues of $Y\in\HD_n$ satisfy $p_0^2=p_1^2$, i.e. that $Y^2=k=p_0^2\in R$. Then the $\HD_n$ modules correspond to 
pairs of Young diagrams. The modules corresponding to $T=(T_0,T_1)$ and $T'=(T_1,T_0)$
are isomorphic.
Furthermore, if $T_0=T_1$ and $2$ is a unit in the ground ring
then the module $V_T$ is $\HD_n$-reducible 
$V_T=V_T^+\oplus V_T^-$. The modules $V_T^+$ and $V_T^-$ 
have the same dimension, but they are not equivalent.
\end{lemma}
\begin{bew}
Consider the map $P:V_T\rightarrow V_{T'}, P(t_0,t_1)=(t_1,t_0)$.
The definition of the action of $X_i,i\geq1$ does obviously
commute with $P$. Furthermore, $i(X_0)=YX_1Y$ acts as $p_0^2X_1$ or $p_1^2X_1$ if
$1,2$ are in the same tableaux 
and as $p_0p_1X_1$ if they are in different tableaux.
Thus  assuming $p_0^2=p_1^2$ we have that $P$ and $i(X_0)$ commute. 
This proves the first claim.
Suppose now that $T_0=T_1$. We see that $P_\pm:=(1\pm P)/2$ 
are projectors on submodules $V_T^\pm$. 
Setting $q=1$ one easily sees that these modules are not isomorphic.
\end{bew}

For the rest of the paper we consider the algebra over the field of fraction of the ground ring.

Since $i$ is an injection and the $V_T$ form a complete set of
simple module of $\HB_n$, we conclude that all simple modules must
be obtained in this way. In fact we have:

\begin{satz}\label{hdmodules}
The set of pairwise non-equivalent absolutely irreducible $\HD_n$
representations is indexed by a set
$I=I_0\cup I_1$, where $I_0:=\{\{T_0,T_1\}\mid T_0\neq T_1,
|T_0|+|T_1|=n\}$
is the set of sets of two Young diagrams with total size $n$ and
$I_1:=\{T^s\mid s\in\{\pm1\},2|T|=n \}$ is the set of $\ZZ_2$ labeled
Young diagrams with exactly $n/2$ boxes. 
\end{satz}
\begin{bew}
$I$ indexes the two types of modules that appeared in the
proof of the preceding lemma.
We show that they are irreducible. This is done by induction
in the same fashion as \cite{ariki}[Proof of 3.9]:
For small $n$ the claim may be shown by direct calculations.
Generators $X_0,\ldots,X_{n-2}$ generate a sub-algebra $\HD_{n-1}'$
of $\HD_n$, which is a surjective image of $\HD_{n-1}$. 
We divide the set of Young tableaux of shape $T$ in two
disjoint sets $T^{(l)},T^{(r)}$ according to the position of $n$.
They split $V_T$ into a direct sum of two  
$\HD_{n-1}'$-submodules and by induction assumption these
are non-equivalent absolutely irreducible representations. 
Now assume that $W\subset V_T$ is a $\HD_n$ submodule.
Viewed as a $\HD_{n-1}'$ module it must contain at least
one of the modules generated by $T^{(l)}$ or $T^{(r)}$.
The action of $X_{n-1}$, however, does not leave such a submodule stable
but allows to generate any other Young tableaux from it. Hence $W$ must
be equal to $V_T$.
The modules are pairwise inequivalent because their restrictions 
to $\HD_{n-1}'$ are by induction assumption inequivalent.
\end{bew}

Note that $I_0$ consists of  half of the $\HB_n$ modules
labeled by two distinct diagrams. Each of the other $\HB_n$ modules
gives rise to two modules of half the dimension. Hence, consistently, we arrive at the
fact  that the dimension of $\HD_n$ is half the dimension of $\HB_n$,
namely $2^{n-1}n!$.

\begin{de}		
Let  $R$ denote an integral domain with units $x,q,\lambda\in R$
such that with  $\delta:=q-q^{-1}$  the relation 
$(1-x)\delta=\lambda-\lambda^{-1}$ holds.
The D-type BMW algebra $\BD_n$
over this ring is the quotient of the $R$-algebra generated by 
$X_0^\pm,\ldots,X_{n-1}^\pm,e_0,\ldots,e_{n-1}$ subject to the following relations:
\begin{eqnarray}
X_iX_jX_i&=&X_jX_iX_j\qquad i,j\geq1,|i-j|=1\label{bd1}\\
X_0X_2X_0&=&X_2X_0X_2\label{bd2}\\
X_0X_1&=&X_1X_0\label{bd3} \\
e_0X_1&=&X_1e_0=e_0\label{bd3a} \\
e_0X_0&=&X_0e_1=e_1\label{bd3b} \\
X_iX_j&=&X_jX_i\qquad |i-j|>1,i,j\geq1\label{bd4} \\
X_ie_i&=&e_iX_i=\lambda e_i\qquad 0\leq i\leq n-1\myownlab{def4}\\
e_iX_j^{\pm1}e_i&=&\lambda^{\mp1}e_i\qquad |i-j|=1\myownlab{def5}\\
e_0X_2^{\pm1}e_0&=&\lambda^{\mp1}e_0\myownlab{def5b}\\
e_2X_0^{\pm1}e_2&=&\lambda^{\mp1}e_2\myownlab{def5v}\\
e_i^2&=&xe_i\qquad 0\leq i\leq n-1\myownlab{lem1a}\\
X_i^{-1}&=&X_i-\delta+\delta e_i\qquad 0\leq i\leq n-1\myownlab{lem1d}\\
e_ie_j&=&e_je_i\qquad |i-j|>1\myownlab{lem1f}\\
e_0e_1&=&e_1e_0\myownlab{lem1fb}\\
e_iX_jX_i&=&X_j^\pm X_i^\pm e_j\qquad |i-j|=1\myownlab{lem1h}\\
e_0X_2X_0&=&X_2^\pm X_0^\pm e_2\myownlab{lem1ha}\\
e_2X_0X_2&=&X_0^\pm X_2^\pm e_0\myownlab{lem1ha}\\
e_ie_je_i&=&e_i\qquad |i-j|=1\myownlab{lem1l}\\
e_0e_2e_0&=&e_0\myownlab{lem1la}\\
e_2e_0e_2&=&e_2\myownlab{lem1lb}
\end{eqnarray}
\end{de}
Most relations are directly transferred from other BMW algebras. And many of them are necessary to render the algebra finite dimensional. An exception is the pair (\ref{bd3a}),(\ref{bd3a}). Even with out them, the algebra is finite dimensional and the structure of the Coxeter graph simply suggests that they should commute (i.e. the first equal sign). However, omitting them would break down the nice relation to the corresponding B-type algebra. The algebra without these relations may, however, be an interesting object for further studies.

\begin{lemma}
The quotient of $\BD_n$ by the ideal generated by $e_0,e_1$ is isomorphic to $\HD_n$.
\end{lemma}
\begin{bew}
The ideal contains all words that contain $e_1$ and thus also $e_2=e_2e_1e_2$ and so forth. 
\end{bew}

The symmetry in the relations shows easily that there is an anti-involution ${}^\ast$ that fixes all generators.

Next, we investigate the relation to the B-type Birman-Murakami-Wenzl algebra $\BB_n$ studied in \cite{ho97}. We assume that we take a common ground ring $R$ for both algebras. For the ground ring of $\BB_n$ it is required to have parameters $q,q_0,q_1,\lambda,x,\delta,A$ with relations $\delta=q-1/q,x\delta=\delta-\lambda+1/\lambda,A\cdot(1-q_0\lambda)=q_1x$. For details see \cite{ho97}. We specialize parameters even further and put a prime on the algebra name to remind about this:

\begin{lemma}
There is an morphism of algebras $i:\BD_n\rightarrow\BB'_n$ to the Birman-Murakami-Wenzl algebra of Coxeter type B with parameters chosen such that $q_0=1/\lambda, q_1=0$ given by $X_0\mapsto \lambda YX_1Y, e_0\mapsto \lambda Ye_1Y, X_i\mapsto X_i, e_i\mapsto e_i, i>0$.
\end{lemma}
\begin{bew}
All calculations are straightforward. We give some examples:

\begin{itemize}
\item $i(X_0X_2X_0)=\lambda^2 YX_1YX_2YX_1Y=\lambda^2 YX_1X_2Y^2X_1Y=\lambda YX_1X_2X_1Y=\lambda  YX_2X_1X_2Y=\lambda  X_2YX_1YX_2=i(X_2X_0X_2)$

\item $i(X_0X_1)=\lambda YX_1YX_1=\lambda X_1YX_1Y=i(X_1X_0)$

\item $i(e_1X_0)=\lambda e_1YX_1Y=e_1=i(e_1)$

\item $i(e_0X_1)=\lambda Ye_1YX_1=Ye_1Y^{-1}=\lambda Ye_1Y=i(e_0)$

\item $i(X_0e_0)=\lambda^2 YX_1YYe_1Y=\lambda YX_1e_1Y=\lambda^2 Ye_1Y=\lambda i(e_0)$

\item The relation $X_i^{-1}=X_i-\delta+\delta e_i$ is tested in its equivalent form $X_i^2=1+\delta X_i-\delta\lambda e_i$:

 $i(X_0^2)=\lambda^2YX_1Y^2X_1Y=\lambda YX_1^2Y=\lambda Y(1+\delta X_i-\delta\lambda e_i)Y=1+\delta i(X_0)-\delta\lambda i(e_0) $

\item $i(e_0X_2e_0)=\lambda^2Ye_1YX_2Ye_1Y=\lambda Ye_1X_2e_1Y= Ye_1Y=\lambda{-1}i(e_0)$

\item $i(e_0^2)=\lambda^2 Ye_1YYe_1Y=\lambda Ye_1e_1Y=x\lambda  Ye_1Y=x\cdot i(e_0)$

\item $0=i(1-X_0^2+\delta X_0-\delta e_0X_0)=1-\lambda^2 YX_1YYX_1Y+\delta \lambda YX_1Y-\delta \lambda^2 Ye_0YYX_0Y=
1-\lambda YX_1X_1Y+\delta \lambda YX_1Y-\delta \lambda Ye_1X_1Y=
1-\lambda Y(1+\delta X_1-\delta \lambda e_1)Y+\delta \lambda YX_1Y-\delta \lambda^2 Ye_1Y
=
1-\lambda Y^2 -\lambda\delta YX_1Y +\delta \lambda^2 Ye_1Y+\delta \lambda YX_1Y-\delta \lambda^2 Ye_1Y
=0$

\item $i(e_0X_1e_0)=Ye_1YX_1Ye_1Y=Ye_1X_1^{-1}X_1YX_1Ye_1Y=\lambda^{-1}Ye_1X_1YX-1Ye_1Y$

$i([e_0,e_1])=0$ follows from $i([X_0,X_1])=0$ and (\ref{bd3}).

\item $i(e_0X_2X_0)=\lambda^2 Ye_1YX_2YX_1Y=\lambda^2 Ye_1Y^2X_2X_1Y=\lambda Ye_1X_2X_1Y=\lambda X_2YX_1Ye_2=i(X_2X_0e_2)$

\item $i(e_0e_2e_0)=\lambda^2 Ye_1Ye_2Ye_1Y=\lambda^2 Ye_1YYe_2e_1Y=\lambda Ye_1e_2e_1Y=\lambda Ye_1Y=i(e_0)$
\end{itemize}
\end{bew}

By composition of $i$ with the Markov trace $\tr:\BB'_n\rightarrow R$ we get a trace on $\BD_n$, i.e. we have the following relations:
\begin{eqnarray}
\tr(1)&=&1\\
\tr(ae_{n-1})&=&x^{-1}\tr(a), a\in\BD_{n-1}\\
\tr(aX_{n-1}^\pm)&=&x^{-1}\lambda^\mp\tr(a)                          \end{eqnarray}

The classical limit ($q=1=\lambda,\delta=0$, $x$ becoming independent of $\lambda$) has a graphical interpretation as Brauer diagrams with dots on it. Each arc cannot have more than one dot (due the quadratic relation of $Y$). The anti-involution considered above amounts to a top-down reflection. The image of $\BD_n$ consists of dotted Brauer diagrams with an even number of dots on it. Vertical arcs in $a$ generate vertical arcs and un-dotted arcs in $aa^\ast$. Horizontal arcs in $a$ on the side joined together produce un-dotted circles in $aa^ast$. In the classical limit the trace is  given 
by closing the strings from the right .

Now we can state an important lemma:
\begin{lemma}
The trace is non-degenerate on   $\BD_n={\rm BP}^k_n$.
\end{lemma}

\begin{bew}
Let $\{v_i\mid i=1,\ldots,k^{n-1}(2n-1)!!\}$ be a linear basis
of dotted Brauer graphs. It suffices to show 
${\rm det}(\tr(v_iv_j^\ast)_{i,j})\neq0$.

The closure of $aa^\ast$ is free of dots.
Now assume that $a$ has  $s$ upper (and hence $s$ lower)
horizontal arcs. Then there are $s$ cycles in $aa^\ast$.
Upon closing, another $s$ cycles are produced from 
the remaining horizontal arcs. The vertical arcs of $a$ form a
permutation and $a^\ast$ contains the inverse permutation. 
Upon closing, these $n-2s$ vertical arcs yield $n-2s$ 
cycles. The closure of   $aa^\ast$ has therefore a total of
 $n$ cycles and  $\tr(aa^\ast)=1$.
For all other dotted Brauer diagrams $b\neq a^\ast$ the closure of $ab$ will have less than $n$ cycles and hence trace $tr(ab)$  is a Laurent polynomial in $x$ that with lower degree in $x$. Hence, there is a unique term with highest power in $x$ and its coefficient can't cancel. Thus the determinant will contain this one term and is therefore not 0. If it is not 0 in the classical limit, it cannot be zero in the general case as well.
\end{bew}

%%%%%%%%%%%%%%%%%%%%%%%%%%%%%%%%%%%%%%%%%%%%%%%%%%%%%%%%%%%%%%%%%%%%%%%%%%%%%%%%%

Next, we determine the structure of  $\BD_n$ in the generic case.
It will turn out to be semi-simple and of dimension $2^{n-1}(2n-1)!!$.

 \begin{satz} \myownlab{hauptsatz}
 Let $R$ be a ring as above and  $K$ its field of fractions.
 The algebra $\BD_n=\BD_n(K)$ is semi-simple and its simple components  are indexed by
the sets $I$ from proposition \ref{hdmodules} for $n, n-2,n-4,...$. 
\end{satz}

The proof uses the same techniques as \cite{we2},\cite{ho97}.

\begin{bew}
Induction starts from $\BD_1$. This algebra has basis $\{1,e_0,X_0\}$ and has three 1-d representations given by the three eigenvalues of $X_0$. 
$\BD_2$ has basis $\{1,e_0,e_1,X_0,X_1,X_0X_1\}$ which is obviously closed under multiplication. $X_0$ and $X_1$ each have 3 eigenvalues so that we could expect to have 9 one-dimensional representations. However, the symmetry between $X_0,X_1$ makes representations depend only on the set of eigenvalues.
 
Assume the proposition is shown by induction for $\BD_n$.

We apply  Jones-Wenzl theory \cite{we1},\cite{we2} to the following inclusion
$\BD_{n-1}\subset\BD_n\subset\BD_{n+1}$. The idempotent is $e=x^{-1}e_n$.
This is possible because $\BD_{n-1},\BD_n$ are semi-simple algebras with a
faithful trace by induction assumption.
All required properties needed for $e$ have already been established.
Jones-Wenzl theory asserts the semi-simplicity of
 the ideal $I_{n+1}$ generated by $e$.
The quotient algebra $\BD_{n+1}/I_{n+1}$ is $\HD_{n+1}$ and is semi-simple. 
Since we work over a field we can conclude
that $\BD_{n+1}$ is semi-simple and that it is isomorphic 
to the direct sum  
$\BD_{n+1}=I_{n+1}\oplus\BD_{n+1}/I_{n+1}=I_{n+1}\oplus\HD_{n+1}$.
Jones-Wenzl theory further implies that the simple components of 
$I_{n+1}$ are indexed by the set of modules of $\HD_{n-1}$.

\end{bew}

\end{document}